\documentclass[reqno]{amsproc}
\usepackage{psfig,latexsym,graphicx, epsfig}
\usepackage{amssymb,amsmath,amsthm,amsopn,pifont} 
\usepackage{subfigure}
\usepackage{multirow}
\newcommand\C{{\mathbb C}}

\newcommand\Z{{\mathbb Z}}
\newcommand\Q{{\mathbb Q}}

\newcommand\QED {$\quad\square$}
\newcommand\mixarr{\;\hbox{\raise 2pt \hbox{$\leftarrow$}}
        \hbox{\lower 1pt \hbox{\hskip -9pt $\to$}}\;}

\def\mod{{\rm mod~}}
\newcommand\Norm{{\rm Norm}}

\newtheorem{theorem}{Theorem}[section]
\newtheorem{lem}[theorem]{Lemma}

\newtheorem{coro}[theorem]{Corollary}

\theoremstyle{definition}

\theoremstyle{remark}
\newtheorem{rem}[theorem]{Remark}
\numberwithin{equation}{section}

\input xy
\xyoption{all}

\begin{document}

\title[Arithmetic of Unicritical Polynomial Maps]{ Arithmetic of Unicritical Polynomial Maps}

\author[John Milnor]{J. Milnor} 
\address{Institute for Mathematical Sciences, Stony Brook University, 
Stony Brook, NY 11794-3660}
\curraddr{}
\email{jack@math.sunysb.edu}
\thanks{I want to thank Thierry Bousch for his help
with this manuscript, and the NSF for its support under grant DMSO757856.    }

\def\IMSmarkvadjust{0 pt}
\def\IMSmarkhadjust{0 pt}
\def\IMSmarkhpadding{0 pt}
\def\IMSpubltext{Published in modified form:}
\def\SBIMSMark#1#2#3{
 \font\SBF=cmss10 at 10 true pt
 \font\SBI=cmssi10 at 10 true pt
 \setbox0=\hbox{\SBF \hbox to \IMSmarkhpadding{\relax}
                Stony Brook IMS Preprint \##1}
 \setbox2=\hbox to \wd0{\hfil \SBI #2}
 \setbox4=\hbox to \wd0{\hfil \SBI #3}
 \setbox6=\hbox to \wd0{\hss
             \vbox{\hsize=\wd0 \parskip=0pt \baselineskip=10 true pt
                   \copy0 \break%
                   \copy2 \break%
                   \copy4 \break}}
 \dimen0=\ht6   \advance\dimen0 by \vsize \advance\dimen0 by 8 true pt
                \advance\dimen0 by -\pagetotal
	        \advance\dimen0 by \IMSmarkvadjust
 \dimen2=\hsize \advance\dimen2 by .25 true in
	        \advance\dimen2 by \IMSmarkhadjust

%
%
  \openin2=publishd.tex
  \ifeof2\setbox0=\hbox to 0pt{}
  \else 
     \setbox0=\hbox to 3.1 true in{
                \vbox to \ht6{\hsize=3 true in \parskip=0pt  \noindent  
                {\SBI \IMSpubltext}\hfil\break
                \input publishd.tex 
                \vfill}}
  \fi
  \closein2
  \ht0=0pt \dp0=0pt
 \ht6=0pt \dp6=0pt
 \setbox8=\vbox to \dimen0{\vfill \hbox to \dimen2{\copy0 \hss \copy6}}
 \ht8=0pt \dp8=0pt \wd8=0pt
 \copy8
 \message{*** Stony Brook IMS Preprint #1, #2. #3 ***}
}

\SBIMSMark{2012/3}{March 2012}{}


\mathsurround = 1pt
\maketitle

This note will study complex polynomial maps of degree $n\ge 2$ with only
one critical point.
Such maps can always be put in the standard normal form
\begin{equation}\label{e1}
	f_c(z)~=~z^n+c 
\end{equation}
by an affine change of coordinate. The connectedness locus, consisting
of all $c$ for which the Julia set of $f_c$ is connected, is sometimes
known as the ``multibrot set''. (Compare \cite{S}.)
 It is not difficult to check that the power $\widehat c=c^{n-1}$
is a complete invariant for the holomorphic conjugacy class of $f_c$.

In \S\ref{s1} we will use the alternate normal form
\begin{equation}\label{e2}
	w~~\mapsto~~g_b(w)~=~(w^n+b)/n~,
\end{equation}
with derivative $g'_b(w)=w^{n-1}$, and use the conjugacy invariant 
$\widehat b=b^{n-1}$. These normal forms 
(\ref{e1}) and (\ref{e2})
 are related by the change of variable formula 
$$w~=~n^{1/(n-1)}z\qquad{\rm with}\qquad
 \widehat b~=~n^n \widehat c\,,$$
and hence $b=n^{n/(n-1)}c$.
(In particular, in the degree two case, $b=\widehat b$ is equal to
 $4c=4\widehat c$.)
\smallskip

If $A$ is any ring contained in the complex numbers $\C$, it will be convenient
to use the non-standard notation $\overline A$ for the \textbf{\textit{ integral
closure}}, the ring consisting of all complex numbers which satisfy a monic 
polynomial equation with coefficients in $A$. (See for example \cite{AM}.)

Section \ref{s1} consists of statements about periodic orbits, which are
 proved in \S\ref{s2}. The last section discusses the
critically finite case. 
\medskip

\section{Periodic Orbits.}\label{s1}
 The  following statement generalizes Bousch \cite{Bo}.
\smallskip

\begin{theorem}\label{t-1.1} 
If $w$ is a periodic
point for the map $g_b$, with multiplier $\mu$, then the rings $\Z[\mu]
\,,~\Z[w]\,,~\Z[b]$, and $\Z[\widehat b]$ all have the same integral closure.
\end{theorem}
\smallskip

 Here are some immediate consequences:
\smallskip

\begin{coro}\label{c-1.2}
If any one of the four numbers $\mu\,,~w\,,~b\,,~\widehat b$
belongs to the ring $\overline\Z$ consisting of all algebraic integers,
then all four of these numbers are algebraic integers. As an example, if 
the map $g_b$ is \textbf{\textit{parabolic}}, that is if the
multiplier of some periodic orbit is a root of unity, then the parameters
$b$ and $\widehat b$ are algebraic integers, hence every periodic point $w$
is an algebraic integer,
and the multiplier $\mu$ of every periodic orbit is an algebraic integer.
\end{coro}
\smallskip

(For a sharper version of this statement, see Remark \ref{r-2.2}.)
\smallskip

\begin{rem}
It would be interesting to understand more generally
which rational maps have the property
that all multipliers are algebraic integers. The family of Latt\`es
maps provides one well known collection of non-polynomial examples.
\end{rem}
\smallskip\smallskip

More generally, if $f:\C\to\C$ is any polynomial map with only one critical
point, then we have the following.

\begin{coro}\label{c- 1.3} 
If $\mu$ and $\mu'$ are the multipliers of two
 periodic orbits for $f$, then the integral closure
$\overline{\Z[\mu]}$ is equal to $\overline{\Z[\mu']}$. It follows that this
integral closure depends only on the holomorphic conjugacy
class of $f$ and not on the particular choice of periodic orbit.
\end{coro}

\begin{figure}[ht!]
\centerline{\psfig{figure=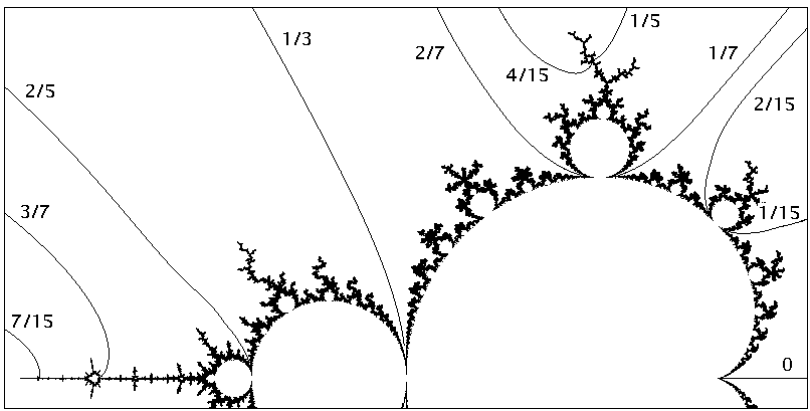,height=2.1in}}
\caption{\label{f1} \sl
Parameter plane for quadratic maps, showing the external rays in the upper
half plane which have period at most four under doubling.}
\end{figure}

\medskip
Now suppose that the parameter value $b$ is the landing
point of an external ray to the connectedness locus
in the $b$-parameter plane, with angle
$p/q\in\Q/\Z$ which is periodic under
multiplication by $n$. (See Figure \ref{f1} for the degree two case.)
Then the associated map
$g_b$ has a parabolic orbit (compare 
\cite{ES}, as well as \cite{DH}, \cite{LS}, \cite{M2}, \cite{S}),
hence the invariant $\widehat b=b^{n-1}=n^n\widehat c$
is an algebraic integer by Corollary \ref{c-1.2}. We will write $b=b(p/q)$, or
$b^{n-1}=\widehat b(p/q)$.

There is a curious relationship between the denominator $q$ of this
angle and the parameter value $b$ or $\widehat b$.
Here are some examples in the quadratic case $n=2$, as illustrated
in Figure \ref{f1}. For the landing points of the $1/3\,,~2/5$, and $3/7$-th 
rays we find that
$$	b(1/3)=-3\,,\qquad b(2/5)=-5\,,\qquad b(3/7)=-7~.$$
At first glance, this relationship between angles and landing points
seems to disappear for the
landing point of the $1/7$-th and $2/7$ rays, with
$~b~=~(-1+3i\sqrt 3)/2~.$
However, this number satisfies the irreducible monic equation
$$	b^2+b+7~=~0~,$$
so the denominator $7=2^3-1$
again appears in the description of the landing point.

Here is a more general statement, working in the $b$ parameter plane for
polynomials of degree $n$. If $b$ has degree $d$ over the rational numbers
$\Q$,  define $\Norm(b)\in\Q$ to be the product of the $d$ algebraic
conjugates of $b$ over $\Q$. Up to sign, this is just the constant term
in the irreducible monic polynomial satisfied by $b$.
If $b$ belongs to the ring $\overline\Z$
of algebraic integers, note that $\Norm(b)\in\Z$.
\smallskip

\begin{theorem}\label{t-1.4} 
 Consider an external ray of angle $p/q$
in this parameter plane with landing point $b=b(p/q)$. If $p/q$ is
periodic under multiplication by $n$ with period $r$, so that $q$ divides
$n^r-1$, and if 
$b$ has degree $d$ over
$\Q$, then it follows that the integer $\Norm(b)$ is a divisor of $(n^r-1)^d$.
Similarly, if $\widehat b$ has degree $\widehat d$ over $\Q$,
 then $\Norm(\widehat b)$ divides $(n^r-1)^{(n-1)\widehat d}\,.$
\end{theorem}
\medskip

Here are some examples:
\smallskip

\begin{itemize}
\item For the landing point of the $1/7$-th ray with ray
period $r=3$ we have $\Norm(b)=b=-7=-(2^3-1)$
\smallskip

\item For the $1/5$-th ray, the ray period is four, and
the irreducible equation is\break $~b^3+9b^2+27b+135=0~$ of degree three
with $|\Norm(b)|=135=3^3\cdot 5$, which is a divisor of $(2^4-1)^3$.
\smallskip

\item For an arbitrary
degree $n\ge 2$, let $b$ be a fixed point of multiplier $\mu$.
Then the equations $g_b'(w)=\mu$ and $g_b(w)=w$
 imply that $w^{n-1}=\mu$ and $b=(n-\mu)w$,
 so that $\widehat b=\mu(n-\mu)^{n-1}$. For $\mu=1$, with $r=1$ and
 $\widehat d=1$, it follows that
 $\widehat b=(n-1)^{n-1}$ is precisely equal to $(n^r-1)^{(n-1)\widehat d}$.
In the case $\mu=-1$, with ray period $r=2$ and $\widehat d=1$,
we get $\widehat b=-(n+1)^{n-1}$, which divides $(n^2-1)^{n-1}$.
\end{itemize}
\smallskip

\section{Proofs.}\label{s2} 
The proofs of the statements of \S\ref{s1} will
depend on some basic properties of the integral
closure. Let $u$ and $v$ be complex numbers. Then clearly
$\overline{\Z[u]}\subset\overline{\Z[v]}$ if and only if 
 $u\in\overline{\Z[v]}$. Note also that
 $\overline{\Z[u^k]}=\overline{\Z[u]}$
for any integer $k>0$. The following statement will be needed.
\smallskip

\begin{lem}\label{l-2.1} 
 Let $u$ and $v$ be complex numbers. If the
  integral closure $\overline{\Z[u]}$ is equal to $\overline{\Z[v]}$,
then it is also equal to $\overline{\Z[uv]}$.
\end{lem}
\smallskip

{\bf Proof.} The product $uv$ certainly belongs to the ring 
$\overline{\Z[u]}\,=\,\overline{\Z[v]}$, hence  $\overline{\Z[uv]}
\subset \overline{\Z[u]}$. Conversely, since $u$ is an element of
$\overline{\Z[v]}$, it satisfies an equation of the form
$$ u^k~=~\sum_{i=0}^{k-1}\sum_{j=0}^\ell n_{i,j}\,u^iv^j\qquad{\rm
  with}\qquad n_{i,j}\in\Z~.$$
 Multiplying both sides of this equation by $u^\ell$, the result
 can be written as
$$ u^{\ell+k}~=~\sum_{i=0}^{k-1}\sum_{j=0}^\ell
n_{i,j}\,u^{\ell+i-j}(uv)^j~,$$
which proves that $u\in\overline{\Z[uv]}$ hence
 $\overline{\Z[u]}\subset\overline{\Z[uv]}$.\QED
\bigskip

{\bf Proof of Theorem~\ref{t-1.1}.}
We can write the $k$-fold iterate $~g_b^{\circ k}(w)~$
as a polynomial with integer coefficients divided by a common
denominator as follows. Set
$$	g_b^{\circ k}(w)~=~ P_k(b,w)/N_k $$
where  $~N_k\,=\,n^{1+n+n^2+\cdots+n^{k-1}}\,=\,n\,N_{k-1}^{\,n}$. Then
$$	P_1(b,w)~=~w^n+b ~,$$
and a straightforward induction shows that
$$	 P_{k+1}(b,w)~=~P_k(b,w)^n\,+\,N_k^{\,n}\, b~.$$
It follows easily that $P_k(b,w)$ is a polynomial  in two variables
with integer coefficients, and that $P_k(b,w)$ is monic of degree $n^k$
when considered as a polynomial in $w$ with coefficients in $\Z[b]$, or
monic of degree $n^{k-1}$ when considered as a polynomial in $b$ with
coefficients in $\Z[w]$.

Now suppose that $w$ is a periodic point for $g_b$, with period $h$.
Then $$g_b^{\circ h}(w)-w=0~,\quad{\rm or~ in~ other~ words}\qquad
	P_h(b,w)\,-\,N_h\,w~=~0~.$$
This last polynomial equation is also monic in either $w$ or $b$, so
it follows that $w\in\overline{\Z[b]}$
and that $b\in\overline{\Z[w]}$. Thus the two rings $\Z[b]$ and $\Z[w]$
have the same integral closure. It follows that the ring
$\Z[\widehat b]$ has this same integral closure.\medskip

Now let 
$$	w=w_0~\mapsto~ w_1~\mapsto~\cdots~\mapsto~ w_h=w_0 $$
be any period $h$ orbit for $g_b$. We know from the argument above
that the rings $\Z[w_j]$ all have the same integral closure. It then
follows inductively from Lemma~\ref{l-2.1} that the ring $\Z[w_1w_2\cdots
w_h]$ has the same integral closure. Since the multiplier of this
orbit can be written as $\mu=(w_1\cdots w_h)^{n-1}$, it follows that
$\Z[\mu]$ also has the same integral closure.\QED\medskip

\begin{rem}\label{r-2.2}
 Here is a supplementary statement. By definition,
an element $w\in\overline\Z$ is relatively \textbf{\textit{prime}} to $n$
if the ideal $w\overline\Z+n\overline\Z$ is equal to $\Z$, or in other
words if $w$ maps to a unit in the quotient ring
$\overline\Z/n\overline\Z$. Now suppose that $w\in\overline\Z$ is
periodic with multiplier $\mu$ under the map $g_b$. If any one of
the four numbers $w,\,\mu,\,b,\,\widehat b$ is prime to $n$, then it
follows that all four of these numbers are prime to $n$. As an
example, if $g_b$ has a parabolic orbit then all of these numbers are
prime to $n$, but if $g_b$ is critically periodic then none of them
is prime to $n$.

To prove this statement, consider an orbit $\{w_j\}$ of period $h$. Then
$nw_{j+1}=w_j^{\,n}+b$, hence $w_j^{\,n}\equiv -b~ ({\rm mod}~n\overline\Z)$.
Taking the product over the orbit elements and raising to the
$(n-1)$-st power, this yields $\mu^n\equiv (-b)^{(n-1)h}$, and the
  conclusion follows easily.
\end{rem}\medskip

{\bf Proof of Theorem~\ref{t-1.4}.}
 Suppose again that $\{w_1\,,\,\ldots\,,\,w_h\}$
is an orbit of period $h$ for $g_b$, with multiplier $\mu=u^{n-1}$ where
$u=w_1w_2\cdots w_h$. Then we have the congruence
$$	nw_{j+1}~=~w_j^{\,n}+b~\equiv~ w_j^{\,n}\quad(\mod b\,\overline\Z)~. $$
In the situation of Corollary~\ref{c-1.2} where $b$ and the $w_j$
belong to the ring $\overline\Z$ of algebraic integers, we can take
the product over $j$ to obtain
$$	n^h\, u~\equiv~u^n\qquad(\mod b\,\overline\Z)~. $$
If $\mu$, and hence $u$, is a unit in the
ring $\overline\Z$, we can divide this congruence by $u$, yielding
\begin{equation}\label{e-4}
	n^h~\equiv~u^{n-1}=\mu\qquad (\mod b\,\overline\Z)~.
\end{equation}
Now suppose that $\mu$ is a primitive $m$-th root of unity. Then raising this
congruence to the $m$-th power, we obtain
$$	n^{hm}~\equiv~1\qquad (\mod b\,\overline\Z)~. $$
Here the product $hm$ is precisely the smallest integer $r$ such that
the iterate $g_b^{\circ r}$ maps each $w_j$ to itself and has derivative $+1$
at each $w_j$.
If an external ray of angle $p/q$ in the $w$-plane lands on $w_j$, then $r$
is precisely equal to the ray period, that is
the period of $p/q$ under multiplication by $n$. (See for example
\cite{M2}.) Using the usual Douady-Hubbard correspondence between parameter 
plane and dynamic plane, at least one of these $p/q$ is also the angle of an
 external ray in the parameter plane which lands on $b$. (Compare \cite{LS}.)
 {\it Thus we see that the ratio $(n^r-1)/b$ is an algebraic integer.}
	
Now taking the product over the $d$ distinct embeddings of the field $\Q(b)$
into $\C$, we see that the rational number
$$	(n^r-1)^d/\Norm(b) $$
belongs to $\overline\Z$, and hence belongs to the ring
$\overline\Z\cap\Q=\Z$. In other words, the integer $\Norm(b)$ is a divisor
of $(n^r-1)^d$, as asserted. A completely analogous argument proves the
corresponding statement for $\Norm(\widehat b)$.\QED\medskip





\section{Postcritically Finite Maps.}\label{s3}
The situation for parameter values corresponding to postcritically finite
maps is rather different. In this case, it is more convenient to work with
the classical normal form of Equation (\ref{e1}), with invariant 
$\widehat c=c^{n-1}= \widehat b/n^n$. The analogue of Corollary \ref{c-1.2}
 in this context is the following.
\smallskip

\begin{lem}\label{l3.1}
If $z$ is periodic under $f_c$, then $z\in\overline\Z$ if and
only if the parameter $c$ or $\widehat c$ belongs to $\overline\Z$. In this
 case, the multiplier $\mu$ of the orbit 
 belongs to the ideal $n^h\overline\Z$, where $h$ is the period.
\end{lem}
\smallskip

The proof is not difficult. (Compare Remark \ref{r-3.2}.)\qed\medskip

In the critically periodic case, it is not hard to show that
$c\in\overline\Z$. This statement can be sharpened as follows for $c\ne 0$.
\smallskip

\begin{theorem}\label{t-3.1} 
 If the orbit of the critical point is
eventually periodic, then $c$ and $\widehat c=c^{n-1}$ are
algebraic integers, with $\Norm(\widehat c)$ dividing $n$. In the special case
where the critical point is actually periodic with period $>1$, we
can sharpen
this statement to say that $\Norm(\widehat c)=\pm 1$.
\end{theorem}
\smallskip

Here are some quadratic examples. If $c=-1$ then the critical point has
period 2, while if $c^3+2c^2+c+1=0$ it has period 3. A number of critically
preperiodic cases are shown in Figure 2, and described further in Tables 
\ref{t1} and \ref{t2}. (Here the transient time is defined to be
the smallest $t$
such that $f_c^{\circ t}(c)$ is periodic.) Note that there is no evident
arithmetic relation between the external angles and the landing point $c$ in
 these cases.
\smallskip

\begin{figure}[ht!]
\centerline{\psfig{figure=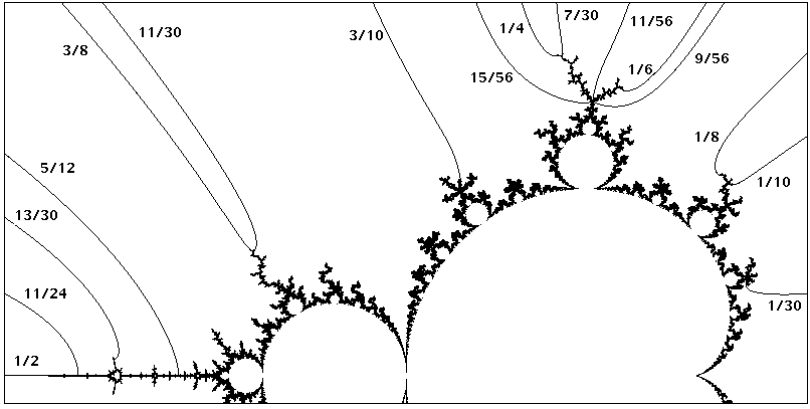,height=2.1in}}\smallskip
\caption{\sl \label{f2} Some eventually periodic rays for the Mandelbrot set.}
\end{figure}

\begin{table}[!ht]
\begin{center}
\begin{tabular}{|ccccc|}
\hline
 $p/q$ & { transient time}~  &~ { eventual period} &~ { degree$(c)$} & 
$|{\rm Norm}(c)|$\cr
\hline
	$\frac{1}{2}$ &  $t= 1$ & $h=1$& $d=1$ &2\cr
	$\frac{1}{6}$  & $t=1$&$h=2$&$d=2$&1\cr
 $\frac{1}{4}, \frac{5}{12}$ & $t=2$&$h=1$&$d=3$&2\cr
$	\frac{1}{8}\,,\,\frac{9}{56}\,,\,\frac{11}{56}\,,\,\frac{15}{56}\,,\,
\frac{ 3}{8}\,,\,\frac{ 11}{24}$ & $t=3$ & $h=1$ & $d=7$ & 2\cr
	$\frac{1}{30}\,,\,\frac{ 1}{10}\,,\,\frac{ 7}{ 30}\,,\,\frac{3}{ 10}\,
,\,\frac{11}{30} \,,\,\frac{13}{ 30}$ &  $t=1$& $h=4$& $d=12$ & 1\cr
\hline
\end{tabular}
\vspace{.3cm}
\caption{\label{t1}\sl Description of the corresponding landing points $c$.}
\end{center}
\end{table}\medskip

\noindent

\begin{table}[!ht]
\begin{center}
\begin{tabular}{|c|}

\hline
	$c+2=0$\cr
	 $c^2+1=0$\cr
	 $c^3+2c^2+2c+2=0$\cr
  $c^7+4c^6+6c^5+6c^4+6c^3+4c^2+2c+2=0$ \cr
  $c^{12}+6c^{11}+15c^{10}+22c^9+23c^8+18c^7+11c^6+8c^5+6c^4+2c^3+1=0\,$. \cr
\hline
\end{tabular}
\vspace{.3cm}
\caption{\sl \label{t2}Corresponding irreducible equations.}
\end{center}

\end{table}
\smallskip

Note that there can be many different postcritically finite parameters
which satisfy the same irreducible equation over $\Q$. This is related to
the fact that the Galois group of $\overline\Q$ over $\Q$ may act in a highly
non-trivial way on these points. (Compare \cite{P}, as well as 
Remark~\ref{r-3.3}.)\bigskip

{\bf Proof of Theorem~\ref{t-3.1}.} Since $c\ne 0$, we can use the
modified normal form
$\zeta\mapsto F_{\widehat c}(\zeta)$ where
$$	F_{\widehat c}(\zeta)~=~\frac{f_c(c\,\zeta)}{c}~=~
\frac{(c\,\zeta)^n+c}{c}~=~ \widehat c\,\zeta^n+1 ~, $$
with critical orbit of the form
$$	0~\mapsto~1~\mapsto~\widehat c+1~\mapsto~\cdots~.$$
The $k$-th point of this
critical orbit can be expressed as a polynomial function $P_k(\widehat c)$,
with $P_1=1$ and
$$P_{k+1}(\widehat c)~=~\widehat c\,P_k(\widehat c)^n+1~.$$
Evidently each
$P_k(\widehat c)$ is a monic polynomial with constant term $P_k(0)=+1$. 
Therefore,
if $F_{\widehat c}$ has periodic critical point, then it follows that $\widehat c$
is a unit in the ring of algebraic integers, with $\Norm(\widehat c)=\pm 1$.

Now suppose that the orbit of zero is eventually periodic but not
periodic. Then the \textbf{\textit{transient time}} $t\ge 1$, 
and the \textbf{\textit{eventual period}} $h\ge 1$ are defined
 as the smallest positive  integers such
that $F_{\widehat c}^{\circ t}(1)$ is periodic of period $h$. It follows that
the two orbit points $P_t(\widehat c)$ and $P_{t+p}(\widehat c)$ are distinct,
and yet have the same image under the $n$-th power map. In other words
the ratio
\begin{equation}\label{e5}
	x~=~P_{t+p}(\widehat c)/P_t(\widehat c) 
\end{equation}
must be an $n$-th root of unity, not equal to $+1$. Hence it must satisfy
the equation
$$	1+x+x^2+\cdots+x^{n-1}~=~ 0~.$$
Clearing denominators, we see that
$$\sum_{i+j=n-1} P_{t+p}(\widehat c)^iP_t(\widehat c)^j~=~ 0~.$$
It is not difficult to check that
 this is a monic polynomial equation in $\widehat c$ with constant
term $n$. Therefore $\widehat c$ is an algebraic integer, and
$\Norm(\widehat c)$ divides $n$.\QED
\smallskip

\begin{rem}\label{r-3.2} 
Note that any periodic point for the map $f_c$ satisfies
a monic equation $f_c^{\circ h}(z)-z=0$ with coefficients in $\Z[c]$. Whenever
$c\in\overline\Z$ and hence $\Z[c]\subset\overline\Z$, it follows 
that $z\in\overline\Z$, hence $f'_c(z)\in n\overline\Z$, so that the multiplier
 $\mu$ belongs to the ideal $n^h\overline\Z$.
In fact there seems to be a strong tendency for periodic points to be units
in the ring $\overline\Z$, so that the ratio $\mu/n^h$ is also a
unit. As an example, suppose that $\widehat
c=-1$ so that the critical orbit has period two. Then the equation
$f_c^{\circ h}(z)-z=0$ is monic with constant term $c=(-1)^{1/(n-1)}$
when $h$ is odd, and the ratio $(f_c^{\circ h}(z)-z)/z$
is monic with constant term $-1$ when $h$ is even. Hence every periodic point
$z\ne 0$ for $f_c$ is an algebraic unit.
\end{rem}\smallskip\smallskip

\begin{rem}\label{r-2.3} Let $\{z_j\}$ be a periodic orbit of period $h>1$
so that
$$ z_{j+1}~=~z_j^{\,n}+c~, $$
where $j$ ranges over $\Z/h\Z$. Using the polynomial expression
$$\phi(x,\,y)~=~\frac{x^n-y^n}{x-y}~=~x^{n-1}+x^{n-2}y+\cdots
+xy^{n-2}+ y^{n-1}~,$$
note the identity
$$ \frac{z_{j+2}-z_{j+1}}{z_{j+1}-z_{j}}
~=~ \frac{z_{j+1}^{\,n}-z_{j}^{\,n}}{z_{j+1}-z_{j}}~=~ \phi(z_j,\,z_{j+1})~.$$
Taking the product over all $j~{\rm modulo}~ h$, it follows that
$~\prod_{j~{\rm mod}~h}\phi(z_j,\,z_{j+1})~=~ 1~.$
(Compare Benedetto \cite{Be}.)
 In particular, if $c\in\overline\Z$ so that the
$z_j$ also belong to $\overline\Z$, then it follows that each expression 
$\phi(z_j,\,z_{j+1})$ is a unit in the ring $\overline\Z$. (For other
``dynamical units'', see \cite{MS}.)
\end{rem}
\medskip

\begin{rem}[Classical Problems]\label{r-3.3}
To conclude this discussion,
we mention two well known unsolved problems. 
\end{rem}
\smallskip

\begin{quote}
{\it If the maps
$f_{c_1}$ and $f_{c_2}$ have parabolic orbits with the same period and the same
ray period, does it follow that the corresponding invariants $\widehat c_1$
and $\widehat c_2$ satisfy the same irreducible equation over $\Q$?\/}
\end{quote}
\smallskip

\noindent In other
words, does it follow that $\widehat c_1$ and $\widehat c_2$ are conjugate under
 some Galois automorphism of the field $\overline\Q$ over $\Q$? 
\smallskip

\begin{quote}
{\it Similarly, if two maps $f_{c_1}$ and $f_{c_2}$ have critical orbits which 
are periodic with the same period, does it follow that $\widehat c_1$ and 
$\widehat c_2$ are Galois conjugate?}
\end{quote}
\smallskip

\noindent
There is a similar question for the eventually periodic case, but the situation
is more complicated. There is an
extra invariant if the degree $n$ is not prime, since the ratio of 
Equation~(\ref{e5}) above
must be a primitive $\tau$-th root of unity for some divisor $\tau$ of $n$,
with $1<\tau\le n$. 
\smallskip

\begin{quote}
{\it If two such parameter values have the same
transient time $t$, the same eventual period $h$, and the same integer
$1<\tau ~|~ n$, does it follow that the corresponding invariants
$\widehat c$ are Galois conjugate?}
\end{quote}

\bigskip

\bibliographystyle{plain}

\begin{thebibliography}{12}


\bibitem[AM]{AM} M. Atiyah and I. Macdonald, ``Introduction to Commutative
Algebra'', Addison-Wesley 1969.

\bibitem[Be]{Be} R. Benedetto, {\it
An elementary product identity in polynomial dynamics,\/}
Amer. Math. Monthly {\bf 108} (2001) 860--864. 

\bibitem[Bo]{Bo} T. Bousch, {\it Les racines des composantes hyperboliques
de $M$ sont des quarts d'entiers alg\'ebriques\/},  (manuscript, 1996). To
appear in ``Frontiers in Complex Dynamics: a volume in honor of John Milnor's 
80th birthday,'' A. Bonifant, M. Lyubich, S. Sutherland, editors.
Princeton University Press 2012.
\smallskip

\bibitem[DH]{DH} A. Douady and J. H. Hubbard, ``\'Etude dynamique des 
polyn\^omes complexes I \& II'', Publ. Math. Orsay (1984-85).\smallskip




\bibitem[LS]{LS} E. Lau and D. Schleicher, {\it Internal addresses in the 
Mandelbrot set and irreducibility of polynomials\/}, Stony Brook IMS Preprint
 1994/19.
\smallskip

\bibitem[ES]{ES} D. Eberlein and D. Schleicher, {\it Rational parameter rays
of multibrot sets\/}, in preparation.

\bibitem[M1]{M1} J. Milnor, ``Dynamics in One Complex Variable'',
Princeton U. Press 2006.\smallskip

\bibitem[M2]{M2} J. Milnor, {\it Periodic Orbits, Externals Rays and the
Mandelbrot Set: An Expository Account\/}, In ``Geometrie Complexe
et Systemes Dynamiques,'' ed.
 M. Flexor, P. Sentenac, J.C. Yoccoz, Ast\'erisque {\bf 261}
(2000) 277--333.\smallskip

\bibitem[MS]{MS} P. Morton and J. Silverman, {\it Periodic points,
 multiplicities, and dynamical units\/}, J. Reine Angew. Math. {\bf 461}
 (1995) 81--122. 

\bibitem[P]{P} K. Pilgrim, {\it Dessins d'enfants and Hubbard Trees\/},
Ann. Sci. \'Ecole Norm. Sup. {\bf (4) 33}  (2000)  671--693.

\bibitem[S]{S} D. Schleicher, {\it
On fibers and local connectivity of Mandelbrot and Multibrot sets\/}, in
 ``Fractal Geometry and Applications: a jubilee of Beno\^\i t
 Mandelbrot. Part 1,''
Proc. Sympos. Pure Math., {\bf 72}, Part 1, Amer. Math. Soc. (2004) 477--517.
\end{thebibliography}

\bigskip\bigskip\bigskip

\rightline{March, 2012}
\end{document}